\documentclass[a4paper,11pt]{amsart}
\usepackage[latin1]{inputenc}
\usepackage{amsmath,amsfonts,amssymb}
\usepackage{color}
\definecolor{gris}{gray}{0.7}

\def\R{\mathbb{R}}
\def\C{\mathbb{C}}
\def\N{\mathbb{N}}

\theoremstyle{plain}
\newtheorem{theo}{Theorem}
\newtheorem{cor}{Corollary}
\newtheorem{prop}{Proposition}
\newtheorem{lm}{Lemma}
\newtheorem{defi}{Definition}

\newenvironment{pr}{\begin{trivlist}\item[]\textit{Proof.}}{\hfill $\square$\end{trivlist}}

\theoremstyle{remark}
\newtheorem{rem}{Remark}

\title[On the Gromov hyperbolicity of the Kobayashi metric]{On the Gromov hyperbolicity of the Kobayashi metric on strictly pseudoconvex regions in the almost complex case}
\author{L\'ea Blanc-Centi}
\address{U.M.P.A., E.N.S. de Lyon\\ 46 all\'ee d'Italie\\ 69364 Lyon cedex 07\\ FRANCE}
\email{Lea.Blanc-Centi@umpa.ens-lyon.fr}

\begin{document}

\begin{abstract}
We prove that every bounded strictly $J$-convex region equipped with the Kobayashi metric is hyperbolic in the sense of Gromov. We apply this result to the study of the dynamics of pseudo-holomorphic maps. 
\end{abstract}

\maketitle

\section{Introduction}
Introduced in 1967 \cite{Kob}, the Kobayashi metric is an important (biholomorphically) invariant metric. It has been used for the study of holomorphic maps and function spaces in several complex variables, but its construction goes through the almost complex situation \cite{KO,Kr}. 
The Kobayashi metric coincides with the Poincar\'e metric on the standard unit disc in $\C$. For the Poincar\'e metric, the unit disc is hyperbolic, and isometries, as well as geodesics, are explicitly known. It is much more difficult to get any information on the global behaviour of the geodesics for the Kobayashi metric on a domain in higher dimension, since there is in general no explicit formula. 

\bigskip

In this paper, we focus on strictly $J$-convex regions in almost complex manifolds, that is, domains whose defining function is strictly $J$-plurisubharmonic. In this case, the local behaviour of the Kobayashi metric has been extensively studyed. Infinitesimally, the Kobayashi metric measures the size of pseudo-holomorphic discs in the domain, and we know various estimates of the infinitesimal Kobayashi metric near the  boundary \cite{Graham,Ma}, even in the almost complex situation \cite{GS,Flo}. Our aim here is to look at the large scale structure of the Kobayashi metric. More precisely, we want to describe the global behaviour of geodesics.

It will be relevant to introduce the notion of $\delta$-hyperbolicity in the sense of Gromov. Let us recall that a geodesic metric space is said to be $\delta$-hyperbolic if the size of every geodesic triangle is less than $\delta$, for some $\delta\ge 0$ only depending on the space. We prove:

\begin{theo}\label{TH2}
Let $J$ be an almost complex structure on $\R^{2n}$ ($n\ge 2$) and $D\subset \R^{2n}$ be a strictly $J$-convex region with connected boundary $\partial D$. Then $D$ equipped with the Kobayashi metric is a geodesic metric space, hyperbolic in the sense of Gromov. Moreover, its boundary as a hyperbolic space is exactly its geometric boundary $\partial D$.
\end{theo}

Notice that the case $n=2$ is due to F. Bertrand \cite{Flo}, and that the link between strict pseudoconvexity and Gromov hyperbolicity was first pointed out for domains in $\C^n$ (that is, when $J$ is the standard complex structure) by Z. Balogh and M. Bonk \cite{BB}. In these two papers, the proof is based on sharp estimates of the Kobayashi metric. In our case, we bypass the lack of such estimates by constructing explicitly a metric $d$ which makes $D$ Gromov hyperbolic (Theorem \ref{TH1}), and for which the geodesics approximate quite roughly ({\em quasi-isometrically}) the geodesics for the Kobayashi metric. The point is that quasi-isometries between {\em geodesic} metric spaces preserve the Gromov hyperbolicity \cite{CDP}. The construction of $d$ uses the contact structure induced on the boundary of the domain by the almost complex structure.


\bigskip

In view of Theorem \ref{TH2}, pseudo-holomorphic maps appear as semi-contractions of Gromov hyperbolic spaces. Taking advantage of this fact, we get (compare with Corollary 20, Chapter 8 of \cite{GdH} which classify the {\em isometries} of Gromov hyperbolic spaces):

\begin{cor}\label{cor:iteration}
Let $J$ be an almost complex structure on $\R^{2n}$ and $D\subset \R^{2n}$ be a strictly $J$-convex region with a connected boundary. Let $F:D\to D$ be a pseudo-holomorphic map. Then either all the orbits (in positive time) of $F$ stay away from the boundary or there is a unique boundary point $p$ such that
$\lim_{k\to +\infty}F^k(x)=p$
for any $x\in D$.
\end{cor}

We stress the fact that this result is obtained only by looking at the geodesics for the Kobayashi metric, that is, by studying the behaviour of pseudo-holomorphic discs. Therefore this method provides results in standard complex analysis without using the algebraic structure of holomorphic functions, but only the geometric properties of the elliptic operator $\bar{\partial}_J$. 

\bigskip

This paper is organized as follows. Section 2 consists of some recalls about the metric notions we need. Section 3 covers the construction of a Gromov hyperbolic metric on any bounded domain in $\R^N$ (see Theorem \ref{TH1}). In Section 4, we focus on the case of strictly $J$-convex domains, and we prove Theorem \ref{TH2} and Corollary \ref{cor:iteration}. 

\bigskip

\smallskip 

\noindent{\em Acknowledgements.}\ 
I would like to thank Bernard Coupet and Herv\'e Gaussier for bringing this question to my attention. I am also grateful to 
\'Etienne Ghys and Herv\'e Pajot, and more generally to 
people who told me about Gromov hyperbolic spaces.

\smallskip

\section{Hyperbolicity in the sense of Gromov}
\subsection{Hyperbolicity}
We begin by giving one of the equivalent definitions for the hyperbolicity in the sense of Gromov:

\begin{defi}\label{def:hyp}
Let $\delta\ge 0$. A metric space $(X,d)$ is {\em $\delta$-hyperbolic} if for every $x,y,z,t\in X$, 
\begin{equation}\label{eqn:hyp}
d(x,y)+d(z,t)\le\mathrm{Max}[d(x,z)+d(y,t), d(x,t)+d(y,z)]+2\delta.
\end{equation}
The metric space $(X,d)$ is said to be {\em hyperbolic in the sense of Gromov} if it is $\delta$-hyperbolic for some $\delta\ge 0$.
\end{defi}

Obviously, every bounded metric space is $\delta$-hyperbolic, for $\delta=\mathrm{diam}X$. The real line is 0-hyperbolic for the Euclidean distance, but the Euclidean space $\R^N$ is not hyperbolic in the sense of Gromov as soon as $N\ge 2$. We refer to \cite{CDP,GdH} for an intensive study of hyperbolic spaces.

For a Gromov hyperbolic space $(X,d)$, one can define a boundary set $\partial_GD$ in the following way (for more details and other constructions, see \cite{GdH}). Fix a point $\omega$ in $X$, and then define the {\em Gromov product} of $x,y\in X$ (with respect to the basepoint $\omega$) as
$$(x,y)_{\omega}=\frac{1}{2}(d(x,0)+d(y,0)-d(x,y)).$$
A sequence $(x_i)$ in $X$ is said to {\em converge at infinity} if $(x_i,x_j)_{\omega}\xrightarrow[{i,j\to\infty}]{}+\infty$. Two sequences $(x_i)$ and $(y_i)$ converging at infinity are {\em equivalent} if $(x_i,y_i)_{\omega}\xrightarrow[{i\to\infty}]{}+\infty$. Theses definitions do not depend on the choice of the basepoint.

\begin{defi}\label{def:GromBd}
The {\em boundary $\partial_G^dX$ of $(X,d)$ as a hyperbolic space} is the set of classes of sequences converging at infinity.
\end{defi}

The Gromov product between two boundary points $a,b\in \partial_G^dX$ is then defined as follows:
\begin{equation}\label{eqn:GromProd}
(a,b)_\omega=\mathrm{Sup}\,\left[\underset{i\to +\infty}{\mathrm{lim\,inf}}(x_i,y_i)_\omega\right]\in(0,+\infty],
\end{equation} 
where the supremum is taken over all sequences $(x_i)$ and $(y_i)$ representing $a$ and $b$, respectively. Then $\partial_G^dX$ carries a canonical topology, for which the balls 
$$\{b\in\partial_G^dX/\ \exp(-(a,b)_\omega)<r\}$$ form a base of neighbourhoods of $a$.

\subsection{Geodesic spaces}
The notion of Gromov hyperbolicity admits a more geometrical characterization in the particular case of a geodesic metric space. Let us recall:

\begin{defi}
A {\em geodesic segment} between two points $x,y$ in a metric space $(X,d)$ is an isometric embedding from $[0;d(x,y)]$ into $X$ connecting $x$ and $y$.

A metric space is called {\em geodesic} if any two points of $X$ can be joined by at least one geodesic segment.
\end{defi}

Even if such a geodesic segment is not necessarily unique, it will be convenient to denote by $\langle x,y\rangle$ the image of $[0;d(x,y)]$ under such an isometry.
In this case, given three points $x,y,z\in X$ there exists a geodesic triangle $\langle x,y\rangle\cup\langle y,z\rangle\cup\langle x,z\rangle$ of vertices $x,y,z$, whose edges are geodesic segments. Such a triangle is said to be {\em $\delta$-thin} if $\mathrm{dist}(w,\langle y,z\rangle\cup\langle z,x\rangle)\le\delta$ for every $w\in\langle x,y\rangle$, and similarly for the other edges. Then we have:

\begin{prop} \cite{GdH}
Let $(X,d)$ be a geodesic metric space. The following conditions are equivalent:
\begin{itemize}
\item $(X,d)$ is hyperbolic in the sense of Gromov;
\item there exists some $\delta\ge 0$ such that every geodesic triangle is $\delta$-thin.
\end{itemize}
\end{prop}

A well-known example of a geodesic space which is Gromov hyperbolic is the unit disc in $\C$ equipped with the Poincar\'e metric (note that in complex dimension one, the Poincar\'e metric is exactly the Kobayashi metric).

\subsection{Length and metric}
We end this preliminary section by giving a criterion for a space to be geodesic. For properties of length spaces, we refer to \cite{GLP}. We begin by defining the length of a {\em path} (that is, a continuous map defined on a segment). 

\begin{defi}\label{def:longueur}
Assume that $(X,d)$ is a metric space. The {\em $d$-length} of a path $\gamma:[a;b]\to X$ is 
\begin{equation}\label{eqn:deflongueur}
\ell(\gamma)=\underset{a=t_0\le\hdots\le t_n=b}{\mathrm{Sup}}\,\sum_{i=1}^nd(t_{i-1},t_i),
\end{equation}
where the supremum (eventually $\infty$) is taken over all possible partitions of $[a;b]$.  
\end{defi}
\noindent Notice that $d(x,y)\le \ell(\gamma)$ for any path $\gamma$ joining $x$ and $y$ in $X$. 

With an additional hypothesis of smoothness, we can give an equivalent definition of the length, which will be useful to make explicit calculations. Let us define the {\em dilation} of $\gamma$ by
$$\mathrm{dil}\,\gamma=\underset{t\not= t'\in[a;b]}{\mathrm{Sup}}\,\frac{d(\gamma(t),\gamma(t'))}{|t-t'|}$$
and the {\em local dilation at point $t$} by $\mathrm{dil}_t\gamma=\underset{\varepsilon\to 0}{\lim}\,\mathrm{dil}\left(\gamma_{|(t-\varepsilon;t+\varepsilon)}\right)$. If we suppose that the path $\gamma$ is Lipschitz, the function $t\mapsto \mathrm{dil}_t\gamma$ is bounded and measurable and 
\begin{equation}\label{eqn:deflongueurbis}
\ell(\gamma)=\int_a^b\mathrm{dil}_t\gamma\,\mathrm{d}t.
\end{equation}
 
\begin{defi}
A metric space $(X,d)$ is a {\em length-space} if the distance between two points of $X$ is the infimum of the lengths of paths joining $x$ and $y$ in $X$.
\end{defi}

\begin{prop}\cite{GLP}\label{prop:compactgeod}
Every compact length-space is geodesic.
\end{prop}

\begin{pr}
It is easy to see that we can consider only paths defined on [0;1] whose parametrization is proportional to the length, that is, 
$$\forall t,t'\in[0;1],\ \ell(\gamma_{|[t;t']})=|t'-t|\,\ell(\gamma).$$
Fix $x,y\in X$ and let $(\gamma_n)_n$ be a sequence of paths in $X$ joining $x$ and $y$ such that $\ell(\gamma_n)\le d(x,y)+1/n$. The family $\{\gamma_n\}_n$ is then equicontinuous and takes its values in $X$. Thus Ascoli's theorem states that, up to extraction, the sequence $(\gamma_n)$ converges uniformly to a path $\gamma$ in $X$ joining $x$ and $y$. Since the sequence $\ell(\gamma_n)$ is uniformly bounded, we get (see for example \cite{Bridson}) 
$$\ell(\gamma)\le\mathrm{lim\,inf}\,\ell(\gamma_n).$$ 
Hence $\ell(\gamma)=d(x,y)$, and $\gamma$ becomes a geodesic segment after reparametrization.
\end{pr}

\section{An example of a Gromov hyperbolic metric on a bounded domain}
Given a metric space $(X,d)$, there always exists a metric that makes it Gromov hyperbolic, namely, $d'=\ln(1+d)$. But this construction leads to a very degenerate situation, since the Gromov boundary $\partial_G^{d'}D$ is either empty or reduced to a single point \cite{CDP}. In fact, the Gromov boundary of a proper hyperbolic metric space can always be obtained as a topological quotient of the geometric boundary \cite{WW}. Here we are interested in the following question: given a bounded domain $D\subset\R^N$, is there a metric on $D$ which makes it Gromov hyperbolic, such that the Gromov boundary $\partial_GD$ coincides with the geometric boundary $\partial D$?
The aim of this section is to give a positive answer.

\begin{theo}\label{TH1} 
Let $D$ be a bounded domain in $\R^{N}$, with $\mathcal{C}^2$-smooth and connected boundary $\partial D$. 
Then there exists a metric on $D$ for which $D$ is geodesic and hyperbolic in the sense of Gromov, and such that the boundary of $D$ as a hyperbolic space coincides (even topologically) with $\partial D$.
\end{theo}

The idea, in order to construct the distance between two points $x,y\in D$, is to give different weights to what happens in the normal direction and in the horizontal (``parallel'' to the boundary) direction. The normal distance will be essentially the Euclidean one, whereas the horizontal distance will be in some way a copy of a metric on the boundary. 

In all this section, we assume that conditions of Theorem \ref{TH1} are satisfied. We will also denote by $d_H$ any metric on the boundary with the following properties:

- the topology induced by $d_H$ on $\partial D$ agrees with the standard topology (hence $(\partial D, d_H)$ is compact);

- the metric space $(\partial D, d_H)$ is geodesic.  

\noindent In view of Proposition \ref{prop:compactgeod}, we can choose for $d_H$ any Riemannian metric, or any sub-Riemannian metric satisfying the conditions of the Chow connectivity theorem (see \cite{Be}, Corollary 2.6; we will precise this point in Section 4).


\subsection{Construction of $d$\label{section:Depsilon}}
In \cite{BS}, Z. M. Balogh and O. Schramm constructed a hyperbolic metric on a one-sided neighbourhood of a bounded metric space $(Y,d_Y)$ (containing more than one point) in the following way. Let $X=Y\times (0;\mathrm{diam}Y]$, then set 
$d_X((p,s),(q,t))=2\ln\left(\frac{d_Y(p,q)+\max(s,t)}{\sqrt{st}}\right)$.
This construction was extended in \cite{BB} to obtain a pseudometric on a $\mathcal{C}^2$-smooth bounded strictly pseudoconvex domain in $\C^n$. 
Here we follow the same idea.

For each $x\in D$, denote by $h(x)$ the squared root of the Euclidean distance of $x$ to the boundary, and choose a point $\pi(x)\in\partial D$ such that $\|x-\pi(x)\|=h(x)^2$.
Then set
\begin{equation}\label{eqn:defg}
\forall x,y\in D,\ g(x,y)=2\ln\left(\frac{d_H(\pi(x),\pi(y))+h(x)\vee h(y)}{\sqrt{h(x)h(y)}}\right)
\end{equation}
where $h(x)\vee h(y)$ denotes the maximum of $h(x)$ and $h(y)$.

There is a certain ambiguity in the definition of $g$ due to the fact that a map $\pi$ with the required properties is not uniquely determined on the whole domain $D$. But since $D$ has a $\mathcal{C}^2$-smooth boundary, one can find some constant $\epsilon>0$ such that for all $p\in\partial D$, the closed ball of radius $\epsilon$ centered at $p-\epsilon \overrightarrow{n}_p$ is in $D\cup\{p\}$ (where $\overrightarrow{n}_p$ denotes the outer normal to $D$ at $p$). Hence, for every point $x\in D$ whose Euclidean distance $h(x)^2$ to the boundary is smaller than $\epsilon$, the projection $\pi(x)$ on $\partial D$ is uniquely defined. Let us fix such an $\epsilon$ for the sequel of the paper.

We set $D_\epsilon=\{x\in D/\ h(x)^2\le\epsilon\}$. Hence the projection $\pi$ is uniquely defined on $D_\epsilon$, and different choices of $\pi$ on $D\setminus D_\epsilon$ lead to functions $g$ that agree near the boundary. We also have that the map $\pi$ is of class $\mathcal{C}^1$ on $D_\epsilon$ (see for example Lemma 2.1 in \cite{BB}).

\begin{figure}
\begin{center}
\begin{picture}(9.2,6)
\qbezier(4,1)(4,1)(7.5,1)
\qbezier(7.5,1)(8.85,1.1)(8.97,3)
\qbezier(7.5,5)(8.85,4.9)(8.97,3)
\qbezier(4,5)(4,5)(7.5,5)
\qbezier(2.2,1.8)(2.92,1.797)(3,1.5)
\qbezier(3,1.5)(3.05,1.05)(4,1)
\qbezier(2.2,4.2)(2.92,4.203)(3,4.5)
\qbezier(3,4.5)(3.05,4.995)(4,5)
\qbezier(0.2,3)(0.27,1.63)(1.5,1.6)
\qbezier(1.5,1.6)(1.69,1.62)(1.8,1.68)
\qbezier(1.8,1.68)(1.95,1.79)(2.2,1.8)
\qbezier(0.2,3)(0.27,4.37)(1.5,4.4)
\qbezier(1.5,4.4)(1.69,4.38)(1.8,4.32)
\qbezier(1.8,4.32)(1.95,4.21)(2.2,4.2)

\qbezier[90](4,2)(6.75,2)(7.5,2)
\qbezier[50](7.5,2)(7.96,2.03)(7.97,3)
\qbezier[50](7.5,4)(7.96,3.97)(7.97,3)
\qbezier[90](4,4)(6.75,4)(7.5,4)
\qbezier[40](2.2,2.8)(3.42,2.785)(3.5,2.4)
\qbezier[20](3.5,2.4)(3.62,2.05)(4,2)
\qbezier[40](2.2,3.2)(3.42,3.215)(3.5,3.6)
\qbezier[20](3.5,3.6)(3.62,3.995)(4,4)
\qbezier[20](1.2,3)(1.21,2.62)(1.5,2.6)
\qbezier[10](1.5,2.6)(1.65,2.61)(1.75,2.68)
\qbezier[20](1.75,2.68)(1.9,2.79)(2.2,2.8)
\qbezier[20](1.2,3)(1.21,3.38)(1.5,3.4)
\qbezier[10](1.5,3.4)(1.65,3.39)(1.75,3.32)
\qbezier[20](1.75,3.32)(1.9,3.21)(2.2,3.2)

\put(5.3,3.1){$_y$}
\put(5.5,2.97){\textbf{.}}
\put(5.25,2.15){$_{y_\epsilon}$}
\put(5.5,1.97){\textbf{.}}
\put(5.3,0.8){$_{\pi(y)}$}
\put(5.5,0.97){\textbf{.}}
\put(5.565,3){\line(0,-1){2}}
\put(5.565,1.15){\line(1,0){0.15}}
\put(5.715,1.15){\line(0,-1){0.15}}

\put(6.6,3.1){$_z$}
\put(6.5,2.97){\textbf{.}}
\put(6.6,4.15){$_{z_\epsilon}$}
\put(6.5,3.97){\textbf{.}}
\put(6.3,5.2){$_{\pi(z)}$}
\put(6.5,4.97){\textbf{.}}
\put(6.565,3){\line(0,1){2}}
\put(6.565,4.85){\line(-1,0){0.15}}
\put(6.415,4.85){\line(0,1){0.15}}

\put(0.7,3.1){$_x$}
\put(0.7,2.97){\textbf{.}}
\put(1.3,3.1){$_{x_\epsilon}$}
\put(1.13,2.97){\textbf{.}}
\put(-0.5,3.1){$_{\pi(x)}$}
\put(0.13,2.97){\textbf{.}}
\put(0.2,3){\line(1,0){1}}
\put(0.35,3.15){\line(-1,0){0.14}}
\put(0.35,3){\line(0,1){0.15}}

\put(3.8,4.5){${D_\epsilon}$}
\put(8.5,3){\vector(-1,0){0.53}}
\put(8.5,3){\vector(1,0){0.47}}
\put(8.4,3.1){$_{\epsilon}$}
\put(8.3,0.8){$\partial D$}

\end{picture}

\figurename{\ 1}
\end{center}
\end{figure}

For every $x\in D$, we set $x_\epsilon=\pi(x)-\epsilon \overrightarrow{n}_{\pi(x)}$, which is the projection of $x$ on the constricted boundary $\{y\in D/\ h(y)^2=\epsilon\}$. Notice that every point of the segment between $x_\epsilon$ and $\pi(x)$ has the same projection $\pi(x)$ on $\partial D$.

\medskip

\begin{rem}\label{rem:ITg}
It is easy to see that for $x,y\in D$, the equality  $g(x,y)=g(x,z)+g(z,y)$ holds if and only if one of the following conditions is satisfied (for simplicity, we assume here that $h(x)\le h(y)$):
\begin{itemize}
\item $h(z)\le h(x)$, and either $g(x,z)=0$ or $g(y,z)=0$;
\item $h(x)\le h(z)\le h(y)$, and either $\pi(z)=\pi(x)$ or $g(y,z)=0$.
\end{itemize}
\noindent Thus, if $x,y,z\in D_\epsilon$ with $z\notin\{x,y\}$, then $g(x,y)=g(x,z)+g(z,y)$ if and only if $h(x)\le h(z)\le h(y)$ and $\pi(z)=\pi(x)$.
\end{rem}

For instance, on Figure 1, $g(x,y)=g(x,x_\epsilon)+g(x_\epsilon,y)$, but $g(x,y)<g(x,y_\epsilon)+g(y_\epsilon,y)$.
Even if the map $g$ restricted to $D_\epsilon\times D_\epsilon$ is a metric, $(D,g)$ is unfortunately not a geodesic metric space because of the previous remark. Nevertheless, $g$ satisfies condition (\ref{eqn:hyp}) of Definition \ref{def:hyp}. Therefore our approach is the following:
\begin{itemize}
\item near the boundary, modify the metric $g$ in order to obtain a geodesic metric $d$; the aim is to construct $d$ such that $-C+g\le d\le g+C$, to preserve condition (\ref{eqn:hyp}). Hence we set: 
$$\forall x,y\in D_\epsilon,\ d(x,y)={\mathrm{Inf}}\,\ell_g(\gamma)$$ 
where the infimum of the $g$-length is taken over all paths joining $x$ and $y$ in $D_\epsilon$. 

\item perturb roughly $g$ inside $D\setminus D_\epsilon$ in order to get $D\setminus D_\epsilon$ as a geodesic metric space; what happens between two points far from the boundary will not thrust in the condition of Gromov hyperbolicity:
$$\forall x,y\in D\setminus D_\epsilon,\ d(x,y)=\begin{array}{|l}\|x-x_\epsilon\|+d(x_\epsilon, y_\epsilon)+\|y-y_\epsilon\|\ \mathrm{if}\ \pi(x)\not=\pi(y)\\
\|x-y\|\ \mathrm{if}\ \pi(x)=\pi(y)\end{array}.$$

\item define $d$ on the whole domain in a logical way for an expected geodesic metric:
$$\forall x\in D_\epsilon,\ y\notin D_\epsilon,\ d(x,y)=d(x,y_\epsilon)+d(y_\epsilon,y)=d(x,y_\epsilon)+\|y-y_\epsilon\|.$$
\end{itemize}
\noindent In particular, $d$ defines a metric on $D_\epsilon$ (since $g$ is a metric on $D_\epsilon$), and
$$\forall x,y\in D_\epsilon,\ g(x,y)\le d(x,y)$$
by definition of the $g$-length. It is not hard to verify, by treating separately the different cases depending on the position of the points, that the inequality $d(x,z)\le d(x,y)+d(y,z)$ is satisfied for every $x,y,z\in D$, and thus

\begin{prop}
The function $d$ gives a metric on $D$.
\end{prop}

Of course, the metric $d$ is {\em not canonically induced} by the choice of the horizontal metric $d_H$ on $\partial D$. It also depends on the choice of the projection $\pi$ and of the parameter $\epsilon$. However, we will precise in Corollary \ref{cor:constructioncanonique} the dependence of our construction on these various choices {\em up to rough-isometry}.

We begin by giving an explicit upper bound for the distance between two points near the boundary. For this, we compute the $g$-length of some special paths, using formula (\ref{eqn:deflongueurbis}).

\begin{lm}\label{lm:segmentvertical}
Let $x,y\in D_\epsilon$ such that $\pi(x)=\pi(y)$. The ``vertical'' path $\gamma^V:[0;d(x,y)]\to D_\epsilon$ defined by $\gamma^V(t)=y+t\,(x-y)/d(x,y)$ verifies 
$$\ell_g(\gamma^V)=\left|\ln\frac{h(x)}{h(y)}\right|.$$
In particular, $\gamma^V$ is a geodesic segment up to reparametrization, and $d(x,y)=g(x,y)$.
\end{lm}

\begin{pr}
Assume for example that $h(y)\le h(x)$. 
Notice that $\forall t,\ \pi(\gamma^V(t))=\pi(\gamma^V(0))$. Hence for all
$t_1<t_2$,\ $g(\gamma^V(t_1),\gamma^V(t_2))=\ln\left(\frac{h(\gamma^V(t_2))}{h(\gamma^V(t_1))}\right)=\frac{1}{2}\ln\left(\frac{h(y)^2+t_2}{h(y)^2+t_1}\right)=\frac{1}{2}\ln\left(1+\frac{t_2-t_1}{h(y)^2+t_1}\right)$
and $\underset{t_2}{\mathrm{Sup}}\,\frac{g(\gamma^V(t_1),\gamma^V(t_2))}{|t_2-t_1|}=\frac{1/2}{h(y)^2+t_1}$. Finally 
$$\mathrm{dil}_t(\gamma^V)=\lim_{\varepsilon\to 0}\frac{1/2}{h(y)^2+(t-\varepsilon)}=\frac{1/2}{h(y)^2+t}$$
and $\ell_g(\gamma^V)=\int_0^{d(x,y)}\mathrm{dil}_t(\gamma)\mathrm{d}t=\ln\left(\frac{h(x)}{h(y)}\right)$. Thus 
$g(x,y)\le d(x,y)\le \ell_g(\gamma^V)=g(x,y)$, which gives the result.
\end{pr}

\begin{lm}\label{lm:segmenthorizontal}
Let $x,y\in D_\epsilon$ such that $h(x)=h(y)$. There exists a path $\gamma^H:[0;d(x,y)]\to D_\epsilon$ at constant height $h(x)^2$ such that
$$\ell_g(\gamma^H)=\frac{2d_H(\pi(x),\pi(y))}{h(x)}.$$ 
We will say that $\gamma^H$ is a {\em short horizontal path} between $x$ and $y$.
\end{lm}

\begin{pr}
Assume that $\alpha:[0;d_H(\pi(x),\pi(y))]$ is a geodesic segment in $(\partial D,d_H)$ joining $\pi(x)$ and $\pi(y)$:
$\forall t,t'\in [0;1],\ d_H(\alpha(t),\alpha(t'))=|t'-t|$. 
We then construct a path in $D$ between $x$ and $y$ by setting 
$$\forall t\in[0;d_H(\pi(x),\pi(y))],\ \gamma^H(t)=\alpha(t)-h(x)^2\overrightarrow{n}_{\alpha(t)}$$
which is the projection of $\alpha$ at height $h(x)^2$.
Thus for all $t_1<t_2$,
$$\frac{g(\gamma^H(t_1),\gamma^H(t_2))}{|t_2-t_1|}=\frac{2}{|t_2-t_1|}\ln\left(\frac{d_H(\alpha(t_1),\alpha(t_2))}{h(x)}+1\right)=\frac{2}{|t_2-t_1|}\ln\left(\frac{|t_2-t_1|}{h(x)}+1\right).$$
This gives $\mathrm{dil}_t(\gamma^H)=\frac{2}{h(x)}$ and hence $\ell_g(\gamma^H)=2d_H(\pi(x),\pi(y))/h(x)$.
\end{pr}

\begin{lm}\label{dadistancebornee}
There exists some constant $C\ge 0$ such that 
\begin{equation*}
\forall x,y\in D,\ -C+g(x,y)\le d(x,y)\le g(x,y)+C.
\end{equation*}
\end{lm}

\begin{pr}
We should consider various cases depending on the relative positions of $x$ and $y$. Notice first that $g$ and $d$ are uniformly bounded on $(D\setminus D_\epsilon)\times(D\setminus D_\epsilon)$, since $\partial D$ is bounded for $d_H$. Moreover, if $x\in D_\epsilon$ and $y\notin D_\epsilon$,
$$d(x,y)-g(x,y)\le d(x,y_\epsilon)+\|y_\epsilon-y\|+g(y_\epsilon,y)-g(x,y_\epsilon)\le d(x,y_\epsilon)-g(x,y_\epsilon) + C$$
$$d(x,y)-g(x,y)\ge d(x,y_\epsilon)-\|y_\epsilon-y\|-g(y_\epsilon,y)-g(x,y_\epsilon)\ge d(x,y_\epsilon)-g(x,y_\epsilon) - C.$$
Hence we only have to study the case $x,y\in D_\epsilon$. Under this hypothesis, we immediately get $g(x,y)\le d(x,y)$ by definition of $d$. It remains to obtain the right part of the inequality.

For simplicity, we assume $h(x)\le h(y)$. As in Figure 2, let $z=\pi(x)-h(y)^2\overrightarrow{n}_{\pi(x)}$: then $h(z)=h(y)$ and $\pi(z)=\pi(x)$.

\begin{figure}
\begin{center}
\begin{picture}(9,5.2)

\qbezier(0.5,3.5)(4.5,4.8)(8.5,3.5)
\put(8.5,3.7){$\partial D$}
\qbezier[90](1.3,1.08)(4.5,2.24)(7.7,1.08)
\put(8.1,2.6){$D_\epsilon$}
{\color{gris} \qbezier(2.15,2.13)(4.5,2.75)(6.85,2.13)}

{\thicklines\put(1.8,3.84){\line(1,-5){0.34}}}
\put(2.07,2.1){\textbf{.}}
\put(2.23,1.35){\textbf{.}}
\put(1.85,2.15){$_{y}$}
{\thicklines\put(7.2,3.84){\line(-1,-5){0.16}}}
{\color{gris} \put(6.73,2.13){\line(1,5){0.18}}}
\put(6.85,3){\textbf{.}}
\put(6.66,2.1){\textbf{.}}
\put(6.515,1.35){\textbf{.}}
\put(6.86,2.15){$_{z}$}
\put(7.08,3){$_{x}$}

\put(2,1.2){$_{y_\epsilon}$}
\put(6.2,1.2){$_{x_\epsilon=z_\epsilon}$}

\put(1.61,3.82){\textbf{.}}
\put(1.35,4.05){$_{\pi(y)}$}
\put(7.01,3.82){\textbf{.}}
\put(6.85,4.05){$_{\pi(x)}$}

\end{picture}
\end{center}
\figurename{\ 2}
\end{figure}

In view of Remark \ref{rem:ITg} and Lemma \ref{lm:segmentvertical}, and since there is a vertical path between $x$ and $z$:
\begin{eqnarray*}
d(x,y)&\le& d(x,z)+d(z,y)\\
 &=&g(x,z)+d(z,y)=(g(x,y)-g(y,z))+d(z,y)=g(x,y)+(d(y,z)-g(y,z)).
\end{eqnarray*}
Hence we only have to prove that $d(y,z)-g(y,z)$ is uniformly bounded for every $y,z\in D_\epsilon$ such that $h(y)=h(z)$. Here we simplify the notations, by setting $d_H=d_H(\pi(y),\pi(z))$ and $h=h(y)$.\\
\textbullet\ First case: $d_H\le h$.

By Lemma \ref{lm:segmenthorizontal}, we get $d(y,z)\le \frac{2d_H}{h}$, and  $d(y,z)-g(y,z)\le 2\frac{d_H}{h}-2\ln\left(1+\frac{d_H}{h}\right)\le 2$.\\
\textbullet\ Second case: $h\le d_H\le\sqrt{\epsilon}$.

Following Figure 3, we consider the path formed with a short horizontal path joining $\pi(y)-d_H^2\overrightarrow{n}_{\pi(y)}$ and $\pi(z)-d_H^2\overrightarrow{n}_{\pi(z)}$, and the two vertical paths $[y;\pi(y)-d_H^2\overrightarrow{n}_{\pi(y)}]$ and $[z;\pi(z)-d_H^2\overrightarrow{n}_{\pi(z)}]$.

\begin{figure}
\begin{center}
\begin{picture}(9,5.2)

\qbezier(0.5,3.5)(4.5,4.8)(8.5,3.5)
\put(8.5,3.7){$\partial D$}
\qbezier[90](1.3,1.08)(4.5,2.24)(7.7,1.08)
\put(8.1,2.6){$D_\epsilon$}
{\color{gris}\qbezier(4.5,2.75)(5.6,2.71)(6.38,2.57)}

{\put(4.5,4.15){\line(0,-1){0.75}}}
{\color{gris}{\put(4.37,3.4){\line(0,-1){0.66}}}}
\put(4.31,3.37){\textbf{.}}
\put(4.43,3.4){$_{y}$}
\put(4.31,2.72){\textbf{.}}
\put(4.31,1.63){\textbf{.}}
\put(4.31,1.5){$_{y_\epsilon}$}
\put(4.31,4.12){\textbf{.}}
\put(4.2,4.35){$_{\pi(y)}$}

{\put(6.5,3.97){\line(-1,-6){0.12}}}
{\color{gris}\put(6.25,3.23){\line(-1,-6){0.105}}}
\put(6.18,3.2){\textbf{.}}
\put(6.32,3.2){$_{z}$}
\put(6.07,2.55){\textbf{.}}
\put(5.89,1.46){\textbf{.}}
\put(5.83,1.33){$_{z_\epsilon}$}
\put(6.305,3.93){\textbf{.}}
\put(6.22,4.13){$_{\pi(z)}$}

\put(4,3.3){\vector(0,1){0.82}}
\put(4,3.3){\vector(0,-1){0.54}}
\put(3.57,3.5){$_{d_H^2}$}

\end{picture}
\end{center}
\figurename{\ 3}
\end{figure}

We get $d(y,z)\le 2\ln\left(d_H/h\right)+2d_H/d_H$ and thus 
$$d(y,z)-g(y,z)\le 2+2\ln\left(\frac{d_H}{h}\right)-2\ln\left(1+\frac{d_H}{h}\right)\le 2.$$
\textbullet\ Third case: $\sqrt{\epsilon}\le d_H$ (and thus $h\le d_H$).

Considering the path formed with the vertical paths $[y;y_\epsilon]$ and $[z;z_\epsilon]$ and a short horizontal path joining $y_\epsilon$ and $z_\epsilon$, we get $d(y,z)\le 2\ln\left(\sqrt{\epsilon}/h\right)+2d_H/\sqrt{\epsilon}$ and
\begin{eqnarray*}
d(y,z)-g(y,z)&\le& \left[2\ln(\sqrt{\epsilon})+2\frac{d_H}{\sqrt{\epsilon}}+2\ln(1/h)\right]-2\ln(1+d_H/h)\\
 &\le&2\ln(\sqrt{\epsilon})+2\frac{d_H}{\sqrt{\epsilon}}-2\ln d_H=2\frac{d_H}{\sqrt{\epsilon}}-2\ln\frac{d_H}{\sqrt{\epsilon}}\le2\frac{M}{\sqrt{\epsilon}}-2\ln\frac{M}{\sqrt{\epsilon}}
\end{eqnarray*}
where the constant $M$ is any uniform upper bound of $d_H$ on the compact set $\partial D\times\partial D$.
This gives the conclusion.
\end{pr}

\medskip

Finally we answer the following question: is this construction canonical, at least in some sense? We need one more definition.

\begin{defi}
Let $d, d'$ be two metrics on some space $X$. We say that $d$ and $d'$ are {\em rough-isometric} if there exists $c\ge 0$ such that
$$\forall x,y\in X, -c+d(x,y)\le d'(x,y)\le d(x,y)+c.$$
\end{defi}

\begin{cor}\label{cor:constructioncanonique}
Assume that $d$ is obtained with the previous construction by choosing a projection $\pi$, a parameter $\epsilon$ and a metric $d_H$ on $\partial D$, and that $d'$ is obtained by making a different choice $(\pi',\epsilon',d_H')$. Then 
$d$ and $d'$ are rough-isometric if and only if $d_H$ and $d_H'$ are bi-Lipschitzly equivalent.
\end{cor}

In particular, if we impose that $d_H$ is any Riemannian metric on $\partial D$, our construction becomes {\em canonical up to rough-isometry}.

\begin{pr}
Let us denote by $g$ and $g'$ the maps obtained by (\ref{eqn:defg}) by choosing respectively $(\pi,\epsilon,d_H)$ and $(\pi',\epsilon',d_H')$. By Lemma \ref{dadistancebornee}, it is equivalent to prove that $d,d'$ are rough-isometric and to prove that there exists $c\ge 0$ such that 
\begin{equation}\label{eqn:onveut}
\forall x,y\in D,\ -c\le g'(x,y)-g(x,y)\le c.
\end{equation}
But $g$ does not depend on the choice of $\epsilon$, and different choices of $\pi$ lead to functions $g$ that agree near the boundary. Hence we just have to get (\ref{eqn:onveut}) in the case when $\epsilon=\epsilon'$ and $\pi=\pi'$. In view of (\ref{eqn:defg}), it reduces to find some constant $c\ge 0$ such that
$$\forall x,y\in D,\ e^{-c/2}\le \frac{d_H'(\pi(x),\pi(y))+h(x)\vee h(y)}{d_H(\pi(x),\pi(y))+h(x)\vee h(y)}\le e^{c/2}.$$
This implies, by considering points $x,y$ that tend to the boundary, that $d_H$ and $d_H'$ are bi-Lipschitzly equivalent. The converse is immediate.
\end{pr}

\subsection{The metric space $(D,d)$}
\begin{prop}
$(D,d)$ is hyperbolic in the sense of Gromov.
\end{prop}

\begin{pr}
This comes from the previous lemma with the same argument as in \cite{BB}. We recall it for the sake of completeness.

Suppose we are given numbers $r_{ij}\ge 0$ such that $r_{ij}=r_{ji}$ and $r_{ij}\le r_{ik}+r_{jk}$ for $i,j,k\in\{1;2;3;4\}$. Then $r_{12}r_{34}\le 4(r_{13}r_{24})\vee(r_{14}r_{23})$. To see this, we may assume that $r_{13}$ is the smallest of the quantities $r_{ij}$ appearing on the right hand side of this inequality. Then $r_{12}\le r_{13}+r_{32}\le 2r_{23}$ and $r_{34}\le r_{31}+r_{14}\le 2r_{14}$. The inequality follows.

Now let $x_i$, $i\in\{1;2;3;4\}$, be four arbitrary points in $\Omega$, and denote by $p_i=\pi(x_i)$ their projections to the boundary and by $h_i$ their heights. Set $d_{ij}=d_H(p_i,p_j)$ and $r_{ij}=d_{ij}+h_i\vee h_j$. Then
$$(d_{1,2}+h_1\vee h_2)(d_{3,4}+h_3\vee h_4)\le 4[(d_{1,3}+h_1\vee h_3)(d_{2,4}+h_2\vee h_4)]\vee[(d_{1,4}+h_1\vee h_4)(d_{2,3}+h_2\vee h_3)],$$
that is, $g(x_1,x_2)+g(x_3,x_4)\le[g(x_1,x_3)+g(x_2,x_4)]\vee[g(x_1,x_4)+g(x_2,x_3)]+2\ln4$.
By Lemma \ref{dadistancebornee}, this gives exactly (\ref{eqn:hyp}) for the metric space $(D,d)$.
\end{pr}

\begin{lm}\label{lm:Depsilongeod}
For all $x,y\in D_\epsilon$, there exists a path $\gamma:[0;1]\to D_\epsilon$ such that $\ell_g(\gamma)=d(x,y)$.
\end{lm}

\begin{pr}
Given any two points $x,y\in D_\epsilon$, we first prove that 
$$\mathrm{Inf}\{\ell_g(\gamma)/\ \gamma\in D_\epsilon\}=\mathrm{Inf}\{\ell_g(\gamma)/\ \gamma\in K_\epsilon(x,y)\},$$ 
where $K_\epsilon(x,y)=\{z\in D/ \mathrm{min}[h(x),h(y)]\le h(z)\le\sqrt{\epsilon}\}$. 

Let $x_0,y_0\in D_\epsilon$ such that $h(x_0)=h(y_0)=h_0$, and assume that $x_1,y_1\in D_\epsilon$ verify $\pi(x_0)=\pi(x_1)$, $h(x_1)\le h_0$ and $\pi(y_0)=\pi(y_1)$, $h(y_1)\le h_0$.
Then 
\begin{eqnarray*}
g(x_1,y_1)=2\ln\left(\frac{d_H(\pi(x_1),\pi(y_1))+h(x_1)\vee h(y_1)}{\sqrt{h(x_1)h(y_1)}}\right)&\ge& 2\ln\left(\frac{d_H(\pi(x_0),\pi(y_0))}{\sqrt{h(x_1)h(y_1)}}+1\right)\\
 &\ge& g(x_0,y_0).
\end{eqnarray*}
Hence if $\gamma$ is a path joining $x_1$ and $y_1$ such that $\forall t,\ h(\gamma(t))\le h_0$, and if $\gamma_0(t)=\pi(\gamma(t))-h_0^2\overrightarrow{n}_{\pi(\gamma(t))}$ is the projection of $\gamma(t)$ at the height $h_0^2$, then 
$$\forall t,t',\ g(\gamma(t),\gamma(t'))\ge g(\gamma_0(t),\gamma_0(t')).$$
By (\ref{eqn:deflongueur}), we obtain $\ell_g(\gamma)\ge \ell_g(\gamma_0)$.

Now let $\gamma$ be any path in $D_\epsilon$ joining $x$ and $y$. Here we suppose for simplicity that $h(x)\le h(y)$. Assume that $h(\gamma(t))<h(x)$ for some $t$. Set $t_1=\mathrm{Inf}\{t> 0/\ h(\gamma(t))<h(x)\}$, and let $t_2>t_1$ be the first time after $t_1$ when the path $\gamma$ is at height $h(x)^2$, as in Figure 4.

\begin{figure}
\begin{center}
\begin{picture}(9,5.2)

\qbezier(0.5,3.5)(4.5,4.8)(8.5,3.5)
\put(0,3.7){$\partial D$}
\qbezier[90](1.3,1.08)(4.5,2.24)(7.7,1.08)
\put(0.1,2.1){$D_\epsilon$}
{\qbezier(1.3,2.8)(4.5,3.84)(7.7,2.8)}

\put(2.04,2.03){\textbf{.}}
\put(1.85,2.14){$_{y}$}
\put(6.83,3){\textbf{.}}
\put(6.85,2.85){$_{x}$}

\put(1.61,3.805){\textbf{.}}
\put(1.35,4.07){$_{\pi(y)}$}
\put(7.01,3.842){\textbf{.}}
\put(6.84,4.07){$_{\pi(x)}$}

\put(7.2,3.4){\vector(1,4){0.1}}
\put(7.2,3.4){\vector(-1,-4){0.1}}
\put(7.32,3.3){$_{h(x)^2}$}

\qbezier(6,3.22)(6.5,2.6)(6.87,3.02)
{\color{gris}\qbezier(4.95,3.32)(5.3,3.9)(6,3.22)}
\qbezier(3.6,3.27)(3.9,1)(4.95,3.31)
{\color{gris}\qbezier(3,3.22)(3.45,4.2)(3.6,3.29)}
\qbezier(2.5,2.8)(2.8,2.6)(3,3.21)
\qbezier(2.1,2.09)(2.2,3)(2.5,2.8)

\put(4.35,2.15){$_{\gamma}$}
\put(5.95,3.17){\textbf{.}}
\put(5.95,3.4){$_{\gamma(t_1)}$}
\put(4.88,3.27){\textbf{.}}
\put(4.3,3.5){$_{\gamma(t_2)}$}

\end{picture}
\end{center}
\figurename{\ 4}
\end{figure}

The previous argument shows that we can replace $\gamma_{|[t_1;t_2]}$ by a short horizontal path joining $\gamma(t_1)$ and $\gamma(t_2)$ at height $h(x)^2$, of $g$-length smaller than $\ell_g(\gamma_{|[t_1;t_2]})$.

By repeating this construction, we finally obtain a new path $\tilde{\gamma}$ in $D_\epsilon$ between $x$ and $y$, such that $\forall t,\ h(\tilde{\gamma}(t))\ge h(x)=\mathrm{min}[h(x),h(y)]$, and $\ell_g(\tilde{\gamma})\le \ell_g(\gamma)$. This immediately proves the claim. 

In view of Proposition \ref{prop:compactgeod}, we just need to show that $K_\epsilon(x,y)$ is compact for the metric $d$. But if $(x_n)$ is a sequence in $K_\epsilon(x,y)$, then $(\pi(x_n))$ is a sequence in the compact set $(\partial D,d_H)$ and $(h(x_n))$ is a sequence in $[\mathrm{min}[h(x),h(y)];\sqrt{\epsilon}\,]$. Up to extraction, we obtain some $p\in\partial D$ and $\lambda\in[\mathrm{min}[h(x),h(y)];\sqrt{\epsilon}\,]$ such that $d_H(p,\pi(x_n))\to 0$ and $h(x_n)\to\lambda$. Set $a=p-\lambda^2\overrightarrow{n}_p$: then $a\in K_\epsilon(x,y)$ and $d_H(\pi(a),\pi(x_n))\to 0$, $h(x_n)\to h(a)$. Hence
$$d(a,x_n)\le\left|\ln\frac{h(a)}{h(x_n)}\right|+2\,\frac{d_H(\pi(a),\pi(x_n))}{h(a)\vee h(x_n)}\xrightarrow[n\to +\infty]{}0,$$
which concludes the proof.
 \end{pr}

\begin{prop}
The metric space $(D,d)$ is geodesic.
\end{prop}

\begin{pr}
According to Definition \ref{def:longueur}, the metric $d$ induces a length function $\ell_d$ on $D$. If some path $\gamma$ takes its values in $D_\epsilon$, then the definition of $d$ implies $\ell_d(\gamma)=\ell_g(\gamma)$ \cite{GLP}. Thus for all $x,y\in D_\epsilon$,  the previous lemma gives a path $\gamma:[0;1]\to D_\epsilon$ such that $\ell_d(\gamma)=d(x,y)$.

Assume that $x,y\in D\setminus D_\epsilon$. If $\pi(x)=\pi(y)$, then $d$ is exactly the Euclidean distance on the vertical segment $\gamma^V$ between $x$ and $y$. Hence $\ell_d(\gamma^V)=\ell_{eucl}(\gamma^V)=\|x-y\|=d(x,y)$. If $\pi(x)\not=\pi(y)$, then  by construction
$$d(x,y)=d(x,x_\epsilon)+d(x_\epsilon,y_\epsilon)+d(y_\epsilon,y).$$
Since $\pi(x)=\pi(x_\epsilon)$, the length of the vertical path between $x$ and $x_\epsilon$ is just $d(x,x_\epsilon)$, and the length of the vertical path between $y$ and $y_\epsilon$ is $d(y,y_\epsilon)$. Hence we only have to find a geodesic path between $x_\epsilon$ and $y_\epsilon$. But this is given by Lemma \ref{lm:Depsilongeod}.
Finally, the path $\gamma$ formed with the two vertical paths $[x;x_\epsilon]$ and $[y;y_\epsilon]$ and a short horizontal path $\gamma_0$ joining $x_\epsilon$ and $y_\epsilon$ satisfies $\ell_d(\gamma)=d(x,y)$.

In the same way, if $x\in D_\epsilon$ and $y\notin D_\epsilon$, then by construction $d(x,y)=d(x,y_\epsilon)+d(y_\epsilon,y)$. Thus the path formed with any geodesic segment joining $x$ and $y_\epsilon$ in $D_\epsilon$ and the vertical path joining $y_\epsilon$ and $y$ is also geodesic.
\end{pr}

\subsection{Boundary of $D$}
We assume for simplicity of notations that $0\in D$ is the basepoint chosen in the definition of the Gromov product.
It derives from Lemma \ref{dadistancebornee} that there exists a constant $C\ge 0$ such that 
\begin{equation}\label{eqn:gromprodG}
\forall x,y\in D,\ -C\le (x,y)_0-\frac{1}{2}(g(x,0)+g(y,0)-g(x,y))\le C.
\end{equation}
Hence
a sequence $(x_i)$ in $(D,d)$ converges at infinity if and only if
$g(x_i,0)+g(x_j,0)-g(x_i,x_j)\xrightarrow[{i,j\to\infty}]{}+\infty$, that is,
 $$\frac{[d_H(\pi(x_i),\pi(0))+h(x_i)\vee h(0)]\times[d_H(\pi(x_j),\pi(0))+h(x_j)\vee h(0)]}{d_H(\pi(x_i),\pi(x_j))+h(x_i)\vee h(x_j)}\xrightarrow[{i,j\to\infty}]{}+\infty.$$
This happens if and only if $(x_i)$ converges with respect to the Euclidean metric to a point in $\partial D$. Each point in $\partial D$ arises as a limit point of a sequence converging at infinity. Moreover, two sequences converging at infinity are equivalent if and only if their limit points on $\partial D$ are the same. 

Assigning to each equivalence class of sequences in $D$ converging at infinity the unique limit point of each sequence in the class, one can identify canonically the Gromov boundary with the geometric boundary as sets.

It remains to compare the topologies of $\partial_G^dD$ and $\partial D$. To prove that they coincide, it suffices to show (using the previous identification) that there exists a constant $C>0$ such that 
$$\forall a,b\in\partial D,\ \frac{1}{C}d_H(a,b)\le\exp(-(a,b)_0)\le Cd_H(a,b).$$ 
Let $a,b\in\partial D$ and let $(x_i)$, $(y_i)$ be sequences in $D$ which converge to $a$ end $b$, respectively. In view of (\ref{eqn:gromprodG}), we obtain that
$$-C'\le\underset{i\to +\infty}{\mathrm{lim\,inf}}(x_i,y_i)_0-\ln\left(\frac{[d_H(a,\pi(0))+h(0)]\times[d_H(b,\pi(0))+h(0)]}{d_H(a,b)}\right)\le C'$$
where $C'=C+\ln(h(0))$. Since the expression $[d_H(a,\pi(0))+h(0)]\times[d_H(b,\pi(0))+h(0)]$ remains uniformly bounded on $\partial D$, we get the conclusion.

\section{Strictly $J$-convex regions equipped with the Kobayashi metric}
\subsection{Almost complex geometry}
First we recall some definitions.

\begin{defi}
An {\em almost complex structure} $J$ on a (real) manifold $M^{2n}$ is a section from $M$ to $End(TM)$, such that $J^2=-Id$.
\end{defi}

\begin{defi}
A map $F:(M,J)\to(M',J')$ of class $\mathcal{C}^1$ between two almost complex manifolds is said to be {\em $(J,J')$-holomorphic} if its differential mapping satisfies $J'\circ dF=dF\circ J$.

If $(M,J)$ is the unit disc in $\C$ (that is, $\Delta\subset\R^2$ equipped with the standard complex structure), we say that $F$ is a {\em $J$-holomorphic disc}. 
\end{defi}

The Kobayashi infinitesimal pseudometric measures the size of pseudo-holomorphic discs:
for every $p\in M$ and every tangent vector $v$ at point $p$, we set
$$K_{(M,J)}(p,v)=\text{inf}\{\alpha>0/\ \exists h\in\mathrm{Hol}^J(\Delta,M)\ \text{s.t.}\ h(0)=p\ \text{and}\ (\partial h/\partial x)(0)=v/\alpha\},$$
which is well-defined according to \cite{NW}.
The function $K_{(M,J)}$ is upper semicontinuous on the tangent bundle of $M$ \cite{Royden,Kr}. We may therefore define the integrated pseudometric $d_K^{(M,J)}$, induced by the infinitesimal pseudometric as follows. Given a curve $\gamma:[a;b]\to M$, $\mathcal{C}^1$-smooth by paths, we define its Kobayashi length as 
$$\ell_K(\gamma)=\int_a^bK_{(M,J)}(\gamma(t),\gamma'(t))\,\mathrm{d}t.$$
This length function induces the {\em Kobayashi pseudometric} by setting, for all $x,y\in M$,
$$d_K^{(M,J)}(x,y)=\text{Inf}\{\ell_K(\gamma)/\ \gamma\ \text{a\ $\mathcal{C}^1$-smooth-by-paths curve in $M$ joining $x$ and $y$}\}.$$

Most of the basic properties of the Kobayashi pseudometric in complex manifolds remains true in the almost complex case, as the decreasing property under the action of pseudo-holomorphic maps:

\begin{prop}\label{prop:distdecr}
Let $F:(M,J)\to(M',J')$ be a $(J,J')$-holomorphic map. Then for every $p,q\in M$,
$$d_K^{(M',J')}(F(p),F(q))\le d_K^{(M,J)}(p,q).$$
\end{prop}

The Kobayashi pseudometric is the maximal pseudometric with this property 
such that it coincides on the standard unit disc in $\C$ with the Poincar\'e metric $d_P$. In fact, $d_K$ is equivalently defined by 
$$d_K^{(M,J)}(x,y)=\mathrm{Inf}\sum_{k=1}^md_P(\zeta_k,\zeta'_k)$$
where the infimum is taken over all chains of pseudo-holomorphic discs $h_k:\Delta\to M$, $k=1,\hdots,m$ satisfying $h_1(\zeta_1)=x,\ h_k(\zeta'_k)=h_{k+1}(\zeta_{k+1})$ and $h_m(\zeta'_m)=y$ \cite{Royden,Kr}.

Unfortunately, $d_K$ defines only a pseudometric in the general case (note that $d_K^\C$ is identically zero). When $d_K^{(D,J)}$ is a metric, the manifold is said to be {\em Kobayashi hyperbolic}. In this case, the topology induced by $d_K$ is the standard one. 

\subsection{Strictly $J$-convex regions}
Under some special conditions of convexity on its boundary, a domain $D$ in an almost complex manifold is Kobayashi hyperbolic.

\begin{defi}
Let $(M,J)$ be an almost complex manifold and $\rho:M\to\R$ a $\mathcal{C}^2$-smooth function. For all $X\in TM$, define $d^c_J\rho(X)=-d\rho(JX)$ and set $\mathcal{L}^J\rho=d(d^c_J\rho)(X,JX)$. The quadratic form $\mathcal{L}^J\rho$ is called the {\em Levi form} of $\rho$.
\end{defi}

\begin{defi}\label{def:strpsc}
A {\em strictly $J$-convex region} in $(M,J)$ is a bounded domain (connected open set) $D=\{\rho<0\}$, where $\rho$ is a $\mathcal{C}^2$-smooth defining function of $D$ whose Levi form is positive definite in a neighbourhood of $D$.
\end{defi}

\noindent This condition does not depend on the choice of the defining function of $D$.

\smallskip

Because of the strict $J$-convexity, the (almost) complex tangential direction and the normal direction do not play the same role. More precisely, for every $p\in\partial D$, the tangent space $T_pM$ splits into 
$T_pM=N_p^J(\partial D)\oplus T_p^J(\partial D)$,
where $N_p^J(\partial D)$ is the $J_p$-invariant subspace spanned by the real normal to $\partial D$ at $p$, and $T_p^J(\partial D)$ is the maximal $J_p$-invariant subspace of $T_p(\partial D)$. Hence, every vector $v\in T_pM$ can be uniquely written as $v=v_N+v_H$, where $v_N$ (resp. $v_H$) is the normal (resp. horizontal) part of $v$. Notice that if $M$ is of (real) dimension $2n$, then $\mathrm{dim}_\R N_p^J(\partial D)=2$ and $\mathrm{dim}_\R T_p^J(\partial D)=2n-2$.

We extend this splitting of the tangent space at some point $x\in D$ in a small neighbourhood of $\partial D$ as follows. Here we restrict ourselves to the case when $D$ is a domain in the Euclidean space (otherwise, the construction would be well-defined only locally), and we use the notations of Subsection \ref{section:Depsilon}. 
For every $t\in(0;\epsilon]$, consider the constricted boundary 
$$S_t=\{p-t\overrightarrow{n}_p/\ p\in\partial D\}$$
formed with the points in $D$ that are at distance $t$ from $\partial D$. 
Every point $x\in D_\epsilon$ is in some $S_t$ (namely, $t=\mathrm{dist}(x,\partial D)$).
Let $N_x^J$ be the $J_x$-invariant subspace spanned by the real normal to $S_t$ at $x$, and $H_x^J$ be the maximal $J_x$-invariant subspace of $T_x(S_t)$. Accordingly, 
$$T_xM=N_x^J\oplus H_x^J$$
and a vector $v\in T_xM$ can be uniquely written as $v=v_N+v_H$, where $v_N\in N_x^J$ and $v_H\in H_x^J$.
We then have precise estimates for the Kobayashi metric:

\begin{theo}\cite{GS,CGS}\label{estimees}
Let $D\subset (R^{2n},J)$ be a strictly $J$-convex region. Then there exist constants $C,C'>0$ such that 
$$\forall x\in D_\epsilon,\ \forall v\in T_xD,\ C\,\left(\frac{\|v_H\|}{h(x)}+\frac{\|v_N\|}{h(x)^2}\right)\le K_{(D,J)}(x,v)\le C'\,\left(\frac{\|v_H\|}{h(x)}+\frac{\|v_N\|}{h(x)^2}\right),$$
where $h(x)=\sqrt{\mathrm{dist}(x,\partial D)}$.
\end{theo}

\begin{rem}
In the complex case \cite{Ma}, and also in the almost complex case if $n=2$ \cite{Flo}, there are in fact much more precise estimates.
\end{rem}

Notice that, by Theorem \ref{estimees}, $d_K^{(D,J)}$ is not only a pseudometric but a metric on $D$. These estimates also give the completeness of the space, using standard integration arguments. Since $(D,d_K)$ is a length space, locally compact by Theorem \ref{estimees}, this implies that it is {\em proper} (that means, the closed balls are compact) and geodesic (\cite{GLP}, Theorem 1.10). Finally, we have obtained:

\begin{prop}
Let $D\subset (\R^{2n},J)$ be a strictly $J$-convex region. Then $D$ equipped with the Kobayashi metric is a proper and geodesic complete metric space.
\end{prop}

\subsection{Proof of Theorem \ref{TH2}}
In this section we assume that conditions of Theorem \ref{TH2} are satisfied.
It remains to show that $(D,d_K)$ is Gromov hyperbolic, and to determine its boundary $\partial_G^{d_K}D$. We first need a definition.

\begin{defi}
A map $\varphi:(X_1,d_1)\to (X_2,d_2)$ between two metric spaces is called {\em quasi-isometric} if there are some constants $C,C'>0$ such that
$$\forall x,y\in X_1,\ -C'+\frac{1}{C}d_1(x,y)\le d_2(\varphi(x),\varphi(y))\le Cd_1(x,y)+C'.$$
\end{defi}

We prove Theorem \ref{TH2} by using the following result, which is a corollary of Theorem 2.2, Chapter 3 in \cite{CDP}.

\begin{theo}\label{theor}\cite{CDP}
Let $X$ be a space equipped with two metrics $d_1$ and $d_2$, such that $(X,d_1)$ and $(X,d_2)$ are geodesic metric spaces. Assume that the identity map between $(X,d_1)$ and $(X,d_2)$ is quasi-isometric. If $(X,d_2)$ is Gromov hyperbolic, then $(X,d_1)$ is also Gromov hyperbolic. Moreover, the boundaries $\partial_G^{d_1}X$ and $\partial_G^{d_2}X$ are canonically homeomorphic.
\end{theo}

The idea is hence to apply Theorem \ref{TH1} to $D$, for a suitable choice of the metric $d_H$ on the boundary. We are looking for $d_H$ verifying that:
\begin{itemize}
\item the topology induced by $d_H$ on $\partial D$ agrees with the standard topology, and the metric space $(\partial D, d_H)$ is geodesic;
\item there exist some constants $C,C'>0$ such that for all $x,y\in D$,
\begin{equation}\label{eqn:qi}
-C'+\frac{1}{C}d_K(x,y)\le 2\ln\left(\frac{d_H(\pi(x),\pi(y))+h(x)\vee h(y)}{\sqrt{h(x)h(y)}}\right)\le Cd_K(x,y)+C'.
\end{equation}
\end{itemize}

The first condition allows us to construct a metric $d$ on $D$ by means of $d_H$ as in Theorem \ref{TH1}. The second condition is exactly the condition of quasi-isometry between $(D,d_K)$ and $(D,d)$ according to Lemma \ref{dadistancebornee} and (\ref{eqn:defg}).

\smallskip

In view of the different roles played in Theorem \ref{estimees} by the complex tangential and the normal directions, we introduce the Carnot-Carath\'eodory metric in the following way.

Recall (see \cite{Bellaiche}) that a {\em contact structure} on the odd-dimensional manifold $\partial D$ is given by a codimension one subbundle $H\subset T(\partial D)$ with non-degenerate curvature form $\Omega:H\wedge H\to T(\partial D)/H$, which can be defined in the following two equivalent fashions.
\begin{itemize}
\item Represent $H$ locally as the kernel of a 1-form, say $\eta$ on $\partial D$, identify $T(\partial D)/H$ with the trivial line bundle and then define $\Omega$ as $d\eta_{|H}$.
\item Define $\Omega(X,Y)$ on pairs of vector fields tangent to $H$ by $\Omega(X,Y)=[X,Y]\,\mathrm{mod}\,H$: it is indeed a 2-form, i.e. $\Omega(aX,bY)=ab\Omega(X,Y)$ for arbitrary smooth functions $a$ and $b$ on $\partial D$.
\end{itemize}

One may think of $H$ as a distinguished set of directions (tangent vectors) in $\partial D$ which we call {\em horizontal}. A 
$\mathcal{C}^1$-smooth curve in $\partial D$ is said to be {\em horizontal} if the tangent vectors to this curve are horizontal. Notice that the horizontal vector fields and their Lie brackets (that is, the commutators of degree less than 2) span the whole tangent space $T_p(\partial D)$ at every point $p$ of $\partial D$. Whence, the Chow connectivity theorem, and its improvement proved by Gromov (see \cite{Bellaiche} p.95 and 120) about the smoothing of horizontal curves, give
 that for every $p,q\in\partial D$, there exists a horizontal $\mathcal{C}^1$-smooth curve in $\partial D$ joining $p$ and $q$. Here we have used the connectivity of the boundary.


\begin{lm}\label{contact}
If $D$ is a strictly $J$-convex region in $(M^{2n},J)$ with $n\ge 2$, the complex tangent bundle $T^J(\partial D)$ is given by a contact structure on $\partial D$.
\end{lm}

\begin{pr}
The assumption on $n$ is equivalent to $T^J(\partial D)\not=\{0\}$.
For $\rho$ being as in Definition \ref{def:strpsc}, we consider the 1-form $\eta=-d^c_J\rho$ on $T(\partial D)$ and set $\Omega=d\eta_{|T^J(\partial D)}$. Then $T^J(\partial D)=\mathrm{Ker}\,\eta$ and for all $X\in T^J(\partial D)$,  
$\Omega(X,JX)=-\mathcal{L}^J\rho(X)$. Assume there exists $X\in T^J(\partial D)$ such that $\Omega(X,Y)=0$ for all $Y\in T^J(\partial D)$. Since the Levi form of $\rho$ is positive definite, we get $X=0$, and  hence $\Omega$ is non-degenerate.  
\end{pr}

\noindent Now we are able to define the associated {\em horizontal} or {\em Carnot-Carath\'eodory metric}:
$$\forall p,q\in\partial D,\ d_H(p,q)={\mathrm{Inf}}\{\mathrm{lengths\ of\ horizontal\ curves\ between\ }p\ \mathrm{and}\ q\}.$$
In particular, $d_H$ is bounded.

\begin{rem}\label{rem:choixdH}
The definition of $d_H$ also involves an auxiliary Riemannian metric on $\partial D$ in order to define the length of a curve. Nevertheless, the choice of this auxiliary metric will not be important for us, since two different Riemannian metrics give bi-Lipschitzly equivalent horizontal metrics, and hence rough-isometric metrics on $D$ (Corollary \ref{cor:constructioncanonique}).
\end{rem} 

The size of balls for $d_H$ can be described quite explicitly. In particular,  
 the topology induced on $\partial D$ by the Carnot-Carath\'eodory metric agrees with the topology induced by any 
Riemannian metric (see for instance \cite{Be}, Corollary 2.6 and Corollary 7.35).
Therefore $\partial D$ is compact also with respect to the horizontal metric, and Proposition \ref{prop:compactgeod} proves that the metric space $(\partial D, d_H)$ is geodesic. 

\begin{rem}
A geodesic segment is not necessarily smooth, so we cannot call it ``horizontal''; nevertheless, its length is well-defined (see Definition \ref{def:longueur}).
\end{rem} 

Thus Theorem \ref{TH1} gives a metric $d$ such that $(D,d)$ is geodesic and Gromov hyperbolic. Moreover, the estimates of Theorem \ref{estimees} give (\ref{eqn:qi}) uniformly on $D$, by the same arguments as in the proof of Corollary 1.3 in \cite{BB}.
Finally, Theorem \ref{theor} leads to the expected result.

Notice that \cite{BB} uses more precise estimates of the Kobayashi metric, due to \cite{Ma} in the complex case, and thus gets the same inequalities with the stronger condition $C=1$. These precise estimates have also been obtained in the almost complex case in \cite{Flo}, but only for the real dimension 4. 

\subsection{Boundary behaviour of proper pseudo-holomorphic maps}
Let $F:D\to D'$ be a proper pseudo-holomorphic map between two regions $D\subset (\R^{2n},J)$ and $D'\subset (\R^{2n},J')$ satisfying conditions of Theorem \ref{TH2}. We proved in \cite{BC} that such a map extends into $\overline{F}:\overline{D}\to\overline{D'}$, in a $\mathcal{C}^1$-smooth way. Moreover  $\overline{F}(\partial D)\subset \partial D'$, and every point $p\in\partial D$ admits a neighbourhood $V$ in $\overline{D}$ such that $\overline{F}_{|V}:V\to\overline{F}(V)$ is a pseudo-biholomorphism. 
Consider the induced boundary map
$$\partial F:\partial D\to \partial D'.$$
It is (at least) $\mathcal{C}^1$-smooth on the compact set $\partial D$, and hence Lipschitz, with respect to any Riemannian metrics. 

Moreover, since $\partial F$ maps horizontal paths in $\partial D$ to horizontal paths in $\partial D'$, we obtain in addition that it is Lipschitz with respect to any Carnot-Carath\'eodory metric $d_H$ and $d'_H$: there exists some constant $C>0$ such that 
$$\forall p,q\in\partial D,\ d'_H(\partial F(p),\partial F(q))\le C d_H(p,q).$$

\subsection{Dynamics of pseudo-holomorphic maps}
Let $J$ be an almost complex structure on $\R^{2n}$ and $D\subset \R^{2n}$ be a strictly $J$-convex region with connected boundary. Let $F:D\to D$ be a pseudo-holomorphic map. By Proposition \ref{prop:distdecr}, $F$ is a {\em semicontraction} of the metric space $(D,d_K)$. Here we are interested in the behaviour of $F$-orbits. 

\begin{defi}
For every $x_0\in D$, the {\em $F$-orbit of $x_0$} is $\overline{\{F^n(x_0)\}}_{n\in\N^*}\subset(D\cup\partial D)$.
\end{defi} 

Since $(D,d_K)$ is proper, conditions of Theorem 11 of \cite{Karlsson} are satisfied. This gives that either every $F$-orbit is bounded in $(D,d_K)$ or every $F$-orbit accumulates only at $\partial D$. More precisely, let us define the {\em limit set} of the $F$-orbit of $x_0$ by setting
$$L^{x_0}(F)=\overline{\{F^n(x_0)\}}_{n\in\N^*}\cap\partial D.$$
The Gromov hyperbolicity of $(D,d_K)$ allows us to apply Proposition 23 of \cite {Karlsson}, and thus we get that in fact $L^{x_0}(F)$ is included in an intersection of singletons depending only of $F$ (the {\em characteristic set} of $F$). This gives all the conclusion of Corollary \ref{cor:iteration}. 
Note that this result was obtained in the complex case in \cite {Karlsson}, without resort to Gromov hyperbolicity, but using some specifically holomorphic arguments.



\address
\email


\end{document}